\documentclass[12pt,leqno]{article}

\usepackage{amscd}
\usepackage{amssymb}
\usepackage[dvips]{graphics}
\usepackage{epic,eepic,epsfig}
\usepackage{amsmath}

\usepackage{amsfonts}
\usepackage[all]{xy}

\def\C{\mathbb{C}}

\def\P{\mathbb{P}}

\def\Q{\mathbb{Q}}

\def\Z{\mathbb{Z}}

\def\A{\mathcal{A}}
\def\B{\mathcal{B}}

\def\F{\mathbb{F}}

\def\L{\mathcal{L}}

\def\HH{\mathcal{H}}

\def\S{\mathcal{S}}

\def\m{\backslash}

\def\proof{\noindent {\bf Proof. } }
\def\endproof{\noindent \hfill $\square$}

\def\ker{\text{\rm Ker }}

\newtheorem{defi}{Definition}[section]

\newtheorem{thm}{Theorem}[section]
\newtheorem{cor}{Corollary}[section]

\newtheorem{remark}{Remark}[section]

\begin{document}

\pagestyle{myheadings}
\markboth{}{}

\title{\bf Hyperplane arrangements, M-tame polynomials and twisted cohomology}

\author{\bf Alexandru Dimca}

\date{}

\maketitle

\section{Introduction}

Let $\A=\{H_1,...,H_d\}$ be an affine essential hyperplane arrangement in
$\C^{n+1}$, see \cite{OT1}, \cite{OT2} for general facts on arrangements.\\
 We set as usual $M=M(\A)=\C^{n+1} \m X$, $X$ being the union  of all the hyperplanes in $\A$. One of the main problems now in hyperplane arrangement theory is to study the cohomology of the complement $M$ with coefficients in some local system $\L$ on $M$, see for instance the introduction and the references in \cite{CDO} as well as \cite{OT2}.\\
A rank one local system $\L$ on $M(\A)$ corresponds to a homomorphism
$$ \pi _1 (M(\A))_{ab}=H_1(M(\A),\Z) \simeq \Z^d \to \C^*$$
i.e. such a local system $\L$ is determined by a collection $\lambda (\L)=
(\lambda_1 (\L),...,\lambda_d (\L))$ of $d$ non-zero complex numbers. Here
$\lambda_j (\L)$ is the monodromy of the local system $\L$ about the hyperplane $H_j$. We call the local system $\L$ equimonodromical if all these monodromies$\lambda_j (\L)$ are the same, i.e. there is $\lambda    \in \C^*$ such that $\lambda_j (\L)= \lambda $  for all $j=1,...,d$. In such a situation we denote the corresponding local system by $\L_{\lambda}$.   \\
We assume in the sequel that $\lambda_j (\L) \ne 1$ for all $j=1,...,d$, the remaining cases being essentially reduced to this one using \cite{C}. Then, there are unique integers $N>1$ and $0<e_j<N$ for $j=1,...,d$ such that
$$ g.c.d(e_1,...,e_d)=1 ~~ and~~ \lambda_j (\L)=exp(2\pi ie_j/N)     $$
 for all $j=1,...,d$. We set $e=(e_1,...,e_d)$.\\
 For any $i=1,...,d$, let $\ell _i=0$ be an equation for the hyperplane $H_i$ and consider the product 
$$f_e=\prod _{i=1,d}\ell_i^{e_i} \in \C[x_0,...,x_n].$$
Let $d_e=e_1+...+e_d$ be the degree of the polynomial $f_e$. When $e_j=1$ for all $j$, then we simply write $f$ for the corresponding product. Note that $deg(f)=d$ and $f=0$ is an equation for the union $X$.\\

When the arrangement $\A$ is central, i.e. $0 \in H_i$ for all $i=1,...,d$, the above polynomial $f$ is homogeneous  and there is a lot of interest in the associated Milnor fiber
$$F=F(\A)=f^{-1}(1)$$
and the corresponding monodromy action $h^q:H^q(F,\C) \to H^q(F,\C)$ coming from the obvious fibration
$$F \to M \to \C^*$$
see for instance \cite{CS}. In particular, it is known that
$$dimH^q(M,\L_{\lambda})=\dim \ker(h^q-\lambda Id)+\dim \ker(h^{q-1}-\lambda Id)$$
see for instance \cite{DN2}. If we denote by $M^*=M^*(\A)$ the quotient $M(\A)/\C^* \subset \P^n$ and if $\lambda ^d=1$, then there is an induced equimonodromical local system $\L_{\lambda}^*$ on  $M^*(\A)$ and we have
$$dimH^q(M^*,\L_{\lambda}^*)=\dim \ker(h^q-\lambda Id)$$
see \cite{CS}.
When the local system $\L$ is not equimonodromical, then one still has an equality
$$dimH^q(M,\L)=\dim \ker(h^q_e-aId)+\dim \ker(h^{q-1}_e-aId)$$
where $a=exp(2\pi i/N)$, $F_e=f_e^{-1}(1)$ and $h_e:F_e \to F_e$ is the corresponding monodromy operator, see  \cite{DN2}.\\

When the  arrangement $\A$ is not central, the usual way to study the cohomology groups $H^*(M(\A),\L)$ is to
identify $\A$ to a projective arrangement $\A_p$ in $\P^{n+1}$ by adding the hyperplane at infinity, hence $|\A_p|=|\A|+1=d+1$, and then study the Milnor fibration of the central arrangement $\B=Cone(\A_p)$ in $\C^{n+2}$ since $M^*(\B)=M(\A)$. This approach has at least two disadvantages:

(i) we have to increase dimensions by one, e.g. if we start with a line arrangement $\A$, the Milnor fiber $F(\B)$ is a surface;

(ii) if we are interested in the special class of equimonodromical local systems $\L_{\lambda}$ and if $a^{d+1}\not =1$, then the local sytem on $M^*(\B)$ naturally associated to  $\L_{\lambda}$ on $M(\A)$, is no longer equimonodromical.\\

The purpose of this note is to introduce a new approach to the study of the affine arrangement $\A$, generalizing the central arrangement case and avoiding the above two problems.
This approach is based on the study of the topology of the function $f:\C^{n+1} \to \C$ and of its monodromy representation, using the tools developed over the years by many authors, see for instance \cite{B},  \cite{NZ1},  \cite{NZ2}, \cite{PZ} and the new progress on Alexander invariants in  \cite{DN2}.

\section {Affine arrangements and M-tame polynomials}

First we recall the notion of an $M$-tame polynomial introduced in \cite{NZ1}
and later studied in \cite{NZ2}, \cite{NS}. For any polynomial $g \in 
 \C[x_0,...,x_n]$ consider the set
$$M(g)=\{x \in \C^{n+1}| grad(g)(x)=c{\overline x} ~~~for~some~c \in \C \}$$
where $grad(g)(x)=(g_0(x),...,g_n(x) ) $, with $g_k$ the partial derivative of $g$ with respect to $x_k$ and ${\overline x}$ is the complex conjugate of $x$.

\begin{defi} We say that the polynomial $g$ is $M$-tame if for any sequence $\{z^k\} \subset M(g)$ with $lim|z^k|=+\infty$ we have $lim|g(z^k)|=+\infty$.
\end{defi}

It is clear that an $M$-tame polynomial can have only isolated singularities
(see also the proof of Corollary \ref{Damon} below). Therefore our polynomial $f$ associated to an affine arrangement cannot be $M$-tame as soon as $n>1$ (except very special cases). Our first result says that this is not a major drawback.

\begin{thm} \label{Mtame}

Let $\A=\{H_1,...,H_d\}$ be an affine essential hyperplane arrangement in
$\C^{n+1}$ given by the equation $f=0$. Then the following hold.

(i) For $n=1$ the polynomial $f$ is $M$-tame.

(ii)  For $n>1$ as well as  for $n=1$ and $d_e >d$, the polynomial $f_e$ is $M_0$-tame in the following sense:
 for any sequence $\{z^k\} \subset M(f_e) \m X$ with $lim|z^k|=+\infty$ we have 
$lim|f_e(z^k)|=+\infty$.
\end{thm}

\proof

The proof of the first claim is easily reduced to the second and we leave it to the reader. The fact that for $n=1$ the polynomial $f$ has a good behaviour at infinity also follows from our discussion in the next section.

The proof of the second claim above is an improved version of the proof of Lemma 4 in \cite{Bo}. Assume that there is a sequence of points
 $\{z^k\} \subset M(f_e) \m X$ with $lim|z^k|=+\infty$ and $lim|f_e(z^k)| \not= +\infty$. Then, by passing to a subsequence, we can assume that  $limf_e(z^k)=b \in \C$.\\
Since the arrangement $\A$ is essential, the set of indices $j$ such that 
$lim|\ell_j(z^k)|=+\infty$ is not empty. By renumbering the hyperplanes, we can assume that $lim|\ell_m(z^k)|=0$ exactly for $1 \leq m \leq q$ with $q \geq 1$
(this set is non-empty since  $limf_e(z^k)=b$). We set 
$$g=\prod _{m \leq q}\ell_m^{e_m} ~~~and~~~~ h=\prod _{m>q}\ell_m^{e_m}.$$
By a linear unitary change of coordinates we can assume that
$$H_1 \cap ...  \cap H_q: x_0=...=x_p=0$$
with $p \leq q-1$. (The unitary requirement is essential, since the condition of $M_0$-tame is a condition of transversality of the fibers of $f_e$ with respect to large spheres centered at the origin, and such spheres being invariant by unitary transformations, it follows that the condition  $M_0$-tame is also invariant.) Then $\ell_m$ for  $1 \leq m \leq q$ is a linear combination of $x_0,...,x_p$ and $g$ is a homogeneous polynomial of degree $q_e=e_1+...+e_q$ in 
$ \C[x_0,...,x_p]$.\\
Now write $z^k=(z_0^k,...,z_n^k)$ in the above fixed coordinate system and hence $z_m^k \to 0$ for  $1 \leq m \leq p$. There is an integer $K>p$ such that 
$|z^k_K| \to +\infty$.\\
Consider the obvious equality
$$\frac {grad(f_e)} {f_e}= \frac {grad(g)} {g}   + \frac {grad(h)} {h} .$$
By passing to a subsequence if necessary, we can assume that $\ell_j(z^k)$ is bounded away from $0$ for $j>q$. It follows that $\frac {grad(h)} {h}    $
 is bounded on the sequence $z^k$.
This implies that for $i>p$, $\frac {f_{e,i}} {f_e}=\frac {h_i} {h}$ 
 is bounded on the sequence $z^k$.

Consider now the equality

$$\frac {\sum_{i=0,p}f_{e,i}(z^k)z^k_i}{f_e(z^k)}=\frac {\sum_{i=0,p}g_i(z^k)z^k_i}{g(z^k)} +\frac {\sum_{i=0,p}h_i(z^k)z^k_i}{h(z^k)}.$$
By Euler formula, the first term in the right hand side is equal to $q_e>0$, while by the above discussion the second term tends to zero. It follows that there is an integer $L \leq p$ such that $|\frac{f_{e,L}(z^k)}{f_e(z^k)}| \to +\infty$.
Since $z^k \in M(f) \m X$ we have
$$|z^k_K||\frac{f_{e,L}(z^k)}{f_e(z^k)}|=|z^k_L||\frac{f_{e,K}(z^k)}{f_e(z^k)}|.$$
This leads to a contradiction, as the left hand side goes to infinity while
 the right hand side goes to zero, by the definition of $K$ and $L$.

\endproof

This result has the following corollaries, saying that essentially $f_e$ behaves like an $M$-tame polynomial. In fact, only the high connectivity of the general fiber $F_e$ of $f_e$ is lost. On the other hand, the defining condition on the multi-index $e$ implies that this general fiber  $F_e$ is connected, see
\cite{DPu}, Remark (I).

\begin{cor} \label{tubes}
For any $t \in \C$ the inclusion of the fiber $F_t =f_e^{-1}(t)$ into the corresponding tube $T_t=f_e^{-1}(D_t)$, with $D_t$ a small disc in $\C$ centered at $t$, is a homotopy equivalence. In particular, both $X=F_0$ and $T_0$ have the homotopy type of a bouquet of $n$-dimensional spheres.

\end{cor}

\proof
The condition of $M$-tame says that the fibers of $f_e$ are transversal to large enough spheres in $\C^{n+1}$ centered at the origin. The weaker condition $M_0$ says the same thing, if we interpret transversality to the special fiber $X=F_0$ in the stratified sense. So the retractions from $T_t$ to $F_t$ obtained in the $M$-tame case in
\cite{NZ1}, \cite{NZ2} by integrating vector fields exist in our case as well.
The fact that $X$ has the homotopy type of a bouquet of $n$-dimensional spheres is well known, see for instance \cite{DP}.

\endproof

Let $\mu(\A)$ denote the number of spheres in the above bouquet. This number is determined by the following result, see \cite{Da} for a different approach.

\begin{cor} \label{Damon}
The function $f_e:M(\A) \to \C$ induced by the polynomial $f_e$ has only isolated singularities and
$$\sum_{x \in M(\A)} \mu(f_e,x)=\mu(\A)=(-1)^{n+1} \chi( M(\A)).$$
\end{cor}

\proof

If $f_e| M(\A)$ would have non-isolated singularities, then we can find $t \in \C^*$ and an irreducible affine algebraic variety $Y \subset F_t \cap Sing(f_e)$
with $\dim Y >0$. Any sequence of points in $Y$ tending to infinity would then contradict the $M_0$-tameness of $f_e$.\\
To complete the proof, we can use the standard trick used already by Broughton in \cite{B} and deduce that $\C^{n+1}$ can be obtained from $T_0$ by adding
$(n+1)$-cells in number equal to the above sum. Then we have just to use the
obvious equalities $\chi(X)=1 +(-1)^n\mu(\A)$ and $\chi(M(\A))=1 -\chi(X)$.

\endproof

The following result explains the interest of this point of view for the
computation of the twisted cohomology of the complement $M(\A)$ with values in a rank one local system $\L$. For basic facts on the monodromy at infinity of polynomials we refer to \cite{DN1}.

\begin{thm} \label{locsys}

(i)  For any integer $k$ such that $0<k<n$, the restriction of the constructible sheaf $R^kf_{e*}\Q$ to $\C^*$ is a local system corresponding to the monodromy operator
$$M^k_e:H^k(F_e,\Q) \to H^k(F_e,\Q).  $$
Here $F_e$ is the general fiber of the polynomial $f_e$ and $M^k_e$ can be taken to be either the monodromy about the fiber $F_0=X$ or, equivalently, the monodromy at infinity
of the polynomial $f_e$.

(ii) Let $\F_e$ by the $\Z$-cyclic covering of $M(\A)$ corresponding to the kernel of the morphism $f_{e*}:\pi_1(M(\A)) \to \pi_1(\C^*)$ and consider $H_n(\F_e,\Q)$ as a
$\Q[t,t^{-1}]$-module in the usual way. Then there is an epimorphism of 
$\Q[t,t^{-1}]$-modules
$$H_n(F_e,\Q) \to H_n(\F_e,\Q)$$
where in the first module the multiplication by $t$ is either the monodromy
 about the fiber $F_0=X$ or the monodromy at infinity
of the polynomial $f$.

\end{thm}

\proof

The first claim follows from the fact that the isolated singularities of $f|M(\A)$ produce  no changes in the topology of the fibers in dimensions $<n$. In particular, the two monodromy operators in the claim (i) above coincide.\\
Using the above construction of $M(\A)$ starting from a punctured tube about $X=F_0$ (which can  also be done starting from a punctured tube about the infinity, i.e. $f^{-1}(\C \m D_R)$, where $D_R$ is a disc in $\C$ of radius $R>>0$
centered at the origin), the proof is similar to the proofs in \cite{DN2}. Easy examples in the case $n=1$ (to be treated in detail in the next two sections) shows that the two monodromy operators in the claim (ii) above do not coincide in general.

\endproof

\begin{cor} \label{Alexinv}

(i) For any integer $k$ such that $0<k<n$,  one has
$$\dim H^k(M(\A),\L)=N(k,a)+N(k-1,a)$$
where $N(k,a)=\dim \ker (M^k_e-aId)$ and $a=exp(2\pi i/N)$.

(ii) $dim H^n(M(\A),\L)\leq N(n,a)+  N(n-1,a)$ and 
$dim H^{n+1}(M(\A),\L)\leq N(n,a)$.

(iii) Both claims (i) and (ii) above hold for the trivial local system $\C_M$
by taking $a=1$.

\end{cor}

\proof
This claim follows from the fact that $M(\A)$ is obtained, exactly as in the proof above, from the punctured tube $T_0^*=T_0 \m X$ by attaching
$(n+1)$-cells, see also \cite{DN2}. It follows that the inclusion $T_0^* \to M(\A)$ induces an isomorphism
$$H^k( M(\A),\L) \simeq H^k(T_0^*,\L)$$
 for $0<k<n$, and hence the result is obtained exactly as the corresponding result for central arrangements mentionned in the Introduction.
For $k=n$ the inequality comes from the epimorphism in Theorem \ref{locsys}, (ii). The last claim is obvious from the previous discussion.

\endproof

\begin{remark}\label{remarks1}  \rm

(i) The $\Q[t,t^{-1}]$-modules $ H_m(\F_e,\Q)$ are exactly the Alexander invariants of the hypersurface $X$ as discussed in \cite{L}, \cite{D2}, \cite{DN2} and, in the case $n=1$, in \cite{K}.

(ii) The $M_0$-tame polynomials have better topological properties than the 
semitame polynomials considered for instance in \cite{PZ}. In particular,
for an  $M_0$-tame polynomial the monodromy at infinity can be realized as the monodromy \`a la Milnor, i.e. the total space can be chosen to be the complement of $X$ in a very large sphere in $\C^{n+1}$ centered at the origin as in the case of $M$-tame polynomials, see \cite{NZ2}.

(iii) It is not clear whether the monodromy operators $M^k_e:H^k(F_e) \to H^k(F_e)  $
 for $0<k<n$ are semisimple. For $k=1$, this is the case for the eigenvalue
$\lambda =1$, see \cite{DS}. 
In the next section we also show that multiplication by $t$ on $H_1(\F,\C)$
is semisimple when $n=1$.
\end{remark}

The following result describes a way to compute the zeta-function
$$Z(f_e,0)(t)=\prod _m(det(Id-tM_{e,0}^m))^{(-1)^m}$$
of the monodromy operator $M_{e,0}$ of the polynomial $f_e$ about the fiber $X=F_0=f_e^{-1}(0)$.

\begin{thm}\label{mono}
The direct image functor $Rf_*$ commutes on the constant sheaf $\C$ to the vanishing cycle functor $\varphi _f$. In particular
$$Z(f_e,0)=\prod _{S \in \S}Z(f_e,x_S)^{\chi (S)}$$
where $\S$ is a constructible regular stratification of $X$ with connected strata such that all the cohomology sheaves $\HH^m(\varphi _f \C)$ are locally constant along the strata of $\S$, $x_S$ is an arbitrary point in the stratum $S$ and
$Z(f_e,x_S)$ is the local zeta-function of the function germ $(f_e,x_S)$.

\end{thm}

\proof

Exactly as in the case of an $M$-tame polynomial treated in \cite{NS},
the direct image functor $Rf_*$ commutes on the constant sheaf $\C$ to the vanishing cycle functor $\varphi _f$. The formula for the zeta-function is similar to the one in the proper case obtained in \cite{GLM} and is treated in detail
for the case of tame polynomials in \cite{D4}.

\endproof

Note that the above commutativity still holds when we replace the functor  $\varphi _f$ by the subfunctor  $\varphi _{f, \lambda}$ which takes only the vanishing cycles corresponding to a fixed eigenvalue $\lambda$. In particular $\varphi _{f, \lambda}\C=0$ implies $N(k,\lambda)=0$ for all $k$. This is an effective way to get vanishing (or upper bound) results for the cohomology groups $H^*(M(\A),\L_a)$, compare to \cite{CDO}, Corollary 16.
In particular, this remark combined with Corollary \ref{Alexinv} yields the following.
\begin{cor} \label{NC}
If $X$ is a normal crossing divisor and $ \lambda_j(\L) \ne 1$ for all $j=1,...,d$, then
$H^q(M(\A),\L)=0$ for all $q < n$.

\end{cor}

\section{Line arrangements (equimonodromical case)}

In this section we assume that $\A$ is an essential line arrangement in the
plane $\C^2$. Let $n_k$ be the number of $k$-fold intersection points in $X$.
The following formulas are easy to deduce.
$$ \chi(M(\A))=\mu(\A)=1-d+\sum_{m \geq 2}n_m(m-1).$$
$$b_1(F)=1-d+\sum_{m \geq 2}n_m(m-1)m.$$
Indeed, the first formula follows from Corollary \ref{Damon} and the additivity
of Euler characteristic with respect to constructible partitions. The second equality comes from the relation
$$b_1(F)=\sum_{x \in \C^2}\mu(f,x)=\chi(M(\A))+\sum_{x \in X}\mu(f,x).$$
Assume that the $d$ lines in $\A$ have $p$ distinct directions and let $k_j$ be the number of lines having the $j$-th direction. A standard computation shows that the genus (of a smooth projective model) of the general fiber $F$ of the defining polynomial $f$ is given by
$$g=genus(F)=\frac{(d-1)(d-2)}{2}-\sum_{j=1,p}\frac{k_j(k_j-1)}{2}.$$
One can determine the resolution graph of $f$ as defined in \cite{ACD}
in a simple way. In fact $X$ intersects the line at infinity $L_{\infty}$ in exactly $p$ points, say $A_1,...,A_p$ (corresponding to the $p$ distinct directions of lines in $X$). Each of these points has to be blown-up, creating thus an exceptional curve $E_j$. The proper transform of $X$ cuts each $E_j$ in exactly $k_j$ points, and each of them has to be blown-up several times to arrive at a dicritic of degree one. Hence the total number of dicritics is
$$\delta(f) = \sum_{j=1,p}k_j=d.$$
This gives the following.

\begin{cor} \label{Kaliman}
Let $n(F_t)$ denote the number of irreducible components of the fiber $F_t$.
Then Kaliman's inequality 
$$\delta(f) -1 \geq \sum _t (n(F_t)-1)$$
is in our situation an equality. In particular, all the fibers $F_t$ for $t \ne 0$ are irreducible.
\end{cor}

It was known that this inequality is an equality when the general fiber $F$ is a rational curve (i.e. $g=0$), see \cite{Ka}, \cite{ACD}, but here we are not in this case
in general, as can easily be verified using the above formula for the genus $g$. One also has $\dim  \ker (M^1_{\infty}-Id) =\delta (f)-1$ for any polynomial $f:\C^2 \to \C$, see \cite{D3}. Therefore the equality $\delta (f)=d$ implies that
$$\dim  \ker (M^1_{\infty}-Id)=b_1(M(\A))-1 .$$  
By Corollary \ref{Alexinv} (iii), we get the same equality when $n>1$.\\

The multiplicity of $f$ along the line  at infinity $L_{\infty}$ is $d$, along the  exceptional curve $E_j$ is $d-k_j$ and then decreases to one for each
 exceptional curve just before a dicritic. Applying A'Campo's formula for the zeta-function as in \cite{ACD} gives the following formula for the
characteristic polynomial of the monodromy at infinity acting on $H^1(F,\C)$.
$$\Delta_{\infty}(t)=(t-1)(t^d-1)^{p-2}\prod_{j=1,p}(t^{d-k_j}-1)^{k_j-1}.$$
Comparing the degree of this polynomial to the previous formula for $b_1(F)$
we get the following relation among the numerical data associated to the line arrangement $\A$.
\begin{cor} \label{numbers}
$$1-d+\sum_{m\geq 2}n_m(m-1)m=(d-1)^2 -\sum _{j=1,p}k_j(k_j-1).$$
\end{cor}
It is also easy to compute the
characteristic polynomial of the monodromy at zero acting on $H^1(F)$. The result is the following.
$$\Delta_{0}(t)=(t-1)^{\mu(\A)}\prod_{m\geq 2}[(t-1)(t^m-1)^{m-2}]^{n_m}.$$
Moreover, in this case the multiplication by $t$ on $H_1(\F,\C)$
is semisimple. Indeed, using Theorem \ref{locsys} we see that the 
 multiplication by $t$ cannot have larger Jordan blocks for the eigenvalue 
$\lambda=1$ since this is the case for the monodromy at infinity, see \cite{D3} and, more generally. \cite{DS}. But the 
 multiplication by $t$ cannot have larger Jordan blocks for the eigenvalue 
$\lambda \not= 1$ since this is the case for the monodromy at zero, all the singularities on $X$ being weighted homogeneous. This proves the final claim in 
 Remark \ref{remarks1} (iii).   \\

Let $\Delta _f$ be the greatest common divisor of the polynomials
$\Delta_{0}$ and $\Delta_{\infty}$. Let $N_f(\lambda)$, 
$N_{0}(\lambda)$ and  respectively $N_{\infty}(\lambda)$ be the multiplicity of
$\lambda$ as a root of the polynomial  $\Delta _f$, $\Delta_{0}$ and 
 respectively  $\Delta_{\infty}$. The following result can be proved exactly as Corollary \ref{Alexinv}.
\begin{cor} \label{upperbound}
For any $\lambda \in \C^*$, $\lambda \ne 1$, we have
$$\dim H^1(M(\A),\L_{\lambda}) \leq N_f(\lambda)=min (N_{0}(\lambda),N_{\infty}(\lambda))         .$$
\end{cor}
It is interesting to compare this upper-bound to the upper-bound obtained in
\cite{CDO}, Theorem 13. Since this latter result applies to equimonodromical rank one local systems on complements of projective line
arrangements in $\P^2$, we have to assume that $\lambda ^{d+1}=1$ such that the
local system $\L_{\lambda}$ is a  equimonodromical local system on the arrangement complement $M(\A_p)$  as explained in the Introduction.
Under this assumtion, it follows that
$$N_{\infty}(\lambda)= \sum_j (k_j-1)$$
where the sum is over all $j$ such that $\lambda ^{k_j+1}=1$. Since $k_j+1$ is exactly the multiplicity of the corresponding projective arrangement  $\A_p$ at the point $A_j$, it follows that $N_{\infty}(\lambda)$ is exactly the upper-bound obtained in \cite{CDO}, Theorem 13  for the arrangement  $\A_p$ and the line at infinity $L_{\infty}$ as a chosen hyperplane.\\
On the other hand, it is easy to see that
$$N_{0}(\lambda)=\sum _m n_m(m-2)$$
where the sum is over all $m \geq 2$ such that $\lambda ^m=1$. The interested reader will have no problem to find explicit examples of line arrangements
showing that both inequalities $N_{\infty}(\lambda)>N_{0}(\lambda)$ and
$N_{0}(\lambda)>N_{\infty}(\lambda)$ are possible. Hence in some cases, the last corollary above gives better upper-bounds that Theorem 13 in 
\cite{CDO} (for any choice of the line at infinity!). One such example (not very interesting) is
$f=xy(x+1)(y+1)(x+y+10)(x+y+11)(x-y+100)(x-y+101)$ and $\lambda$ a cubic root of unity. Here any line in the associated projective arrangement contains at least a triple point (and hence $N_{\infty}(\lambda) \geq 1$ for any choice of the line at infinity)  , but $X$ has only normal crossings and hence
$N_{0}(\lambda)=0$.

\section{Line arrangements (general case)}

In this section we continue to use the notation from the previous section, in particular $X \cap L_{\infty} =\{A_1,...,A_p\}$. These $p$ line directions induce a partition $(I_1,...,I_p)$ of the set of indices $\{1,...,d\}$ such that
$i \in I_j$ if and only if $H_i \cap  L_{\infty}=A_j$.
Let $C_t=\overline F _t $ be the closure in $\P^2$ of the fiber $F_t=f_e^{-1}(t)$.
Then $C_t$ has exactly $p$ singularities along the line at infinity (namely at the points $\{A_1,...,A_p\}$), and
an easy computation using the additivity of Milnor numbers under a blow-up,
see \cite{D1}, Proposition (10.27) shows that
$$ \mu(C_t,A_j)=d_e(d_j-k_j)+d_j(k_j-2) +1.$$
Here $d_j=\sum_{i \in I_j}e_i$ and $k_j=|I_j|$. This formula implies in the usual way the following equality
$$b_1(F_e)=1+d_e(d-1)-\sum_{j=1,p}d_jk_j.$$
One surprizing consequence of this formula when compared to Corrolary \ref{Damon} is that for a fixed arrangement $\A$ we have $sup_e b_1(F_e)=\infty$, i.e. the topology of the the general fiber $F_e$ becomes more and more complicated as the multiplicities $e$ increase.

Similar considerations as in the previous section shows that $\delta(f_e)=d$, hence the Kaliman's inequality is an equality in this case as well and all the fibers $F_t=f_e^{-1}(t)$ are irreducible for $t \ne 0$. Moreover, we get the following formula for the characteristic polynomial of the monodromy operator $M_{e,\infty}^1$ at infinity of the polynomial $f_e$.
$$\Delta _{e,\infty}(t)=(t-1)(t^{d_e}-1)^{p-2}\prod _{j=1,p}(t^{d_e-d_j}-1)^{k_j-1}.$$
Moreover, Theorem \ref{mono} can be applied in this situation and yields the following formula for the characteristic polynomial of the monodromy operator $M_{e,0}^1$ about the fiber $F_0=X$ of the polynomial $f_e$.
$$\Delta _{e,0}(t)=(t-1)\prod_{lines} (t^{e_j}-1)^{-\chi(H_j^0)}
\prod_{vertices}(t^{d(I_v)}-1)^{|I_v|-2}$$
where the first product is over all the lines $H_j$ and $H_j^0=H_j \m \cup_{i \ne j}H_i$ and the second product is over all the vertices $v$, $I_v$ denotes the set of $m$ such that $v \in H_m$ and $d(J)=\sum_{m \in J}e_m$.

Let us investigate the multiplicity of a root $a=exp(2\pi i/N)$ in these two polynomials $\Delta _{e,\infty}$ and $\Delta _{e,0}$ under the assumption that
$\lambda_j(\L)=a^{e_j} \ne 1$ for any $j$.
Using the above formula for  $\Delta _{e,0}$ it is easy to see that this multiplicity is
$$N_0(a)=mult(a,\Delta _{e,0})=\sum_{vertices}(|I_v|-2)$$
where the sum is over all vertices $v$ in $\C^2$ such that $\prod _{j \in I_v}
\lambda_j(\L)=1$. In a similar way
$$N_{\infty}(a)=mult(a,\Delta _{e,\infty})=\sum_{vertices}(|I_v|-2)$$
where the sum is over all vertices $v \in L_{\infty}$ of the corresponding projective arrangement $\A_p$ in $\P^2$ such that $\prod _{j \in I_v}
\lambda_j(\L_p)=1$, $\L_p$ being the local system $\L$ regarded as a local system on $M(\A_p)$.
Then we have the following result.

\begin{cor}
With the above notation, for any rank one local system $\L$ on $M(\A)$ such that $\lambda _j(\L) \ne 1$ for all $j$ one has
$$\dim H^1(M(\A),\L) \leq min(N_0(a),N_{\infty}(a)).$$
\end{cor}

 The upper-bound on $\dim H^1(M(\A),\L)$ obtained from $N_{\infty}(a)$ can be considered as a generalization of Theorem 13 in \cite{CDO}, which applies only to equimonodromical local systems.

On the other hand, it is easy to give a sheaf theoretic proof of the above Corollary. Indeed, the setting in the proof of Theorem 13 in \cite{CDO}
gives by a slight modification the  upper-bound obtained from $N_{\infty}(a)  $. To get the  upper-bound $N_0(a)$, it is enough to play the same game of comparing the direct image $Rj_*\L$ with the direct image with compact supports 
$Rj_!\L$ as in
 \cite{CDO}, but replacing the affine space $\C^2$ by a large closed ball
$B$ centered at the origin of $\C^2$ and taking $j$ to be the inclusion
$M(\A) \cap B \to B$. Indeed, it is known that the inclusion
$M(\A) \cap B \to M(\A)$ is a homotopy equivalence, see for instance \cite{D2}  p. 26 and hence
$H^1( M(\A),\L) \simeq H^1( M(\A) \cap B,\L)$. Further details will be given elsewhere.
\bigskip

\bibliographystyle{amsalpha}

\bigskip

Laboratoire d'Analyse et G\'eom\'etrie,\\

Universit\'e Bordeaux I,\\

33405 Talence Cedex, FRANCE\\

email: dimca@math.u-bordeaux.fr

\end{document}